\documentclass[11pt,leqno]{amsart} 
\usepackage{graphics}
\usepackage{hyperref}

\usepackage{lipsum}

\let\OLDthebibliography\thebibliography
\renewcommand\thebibliography[1]{
  \OLDthebibliography{#1}
  \setlength{\parskip}{2pt}
  \setlength{\itemsep}{2pt plus 0.3ex}
}

\newtheorem{thm}{Theorem}[section]
\newtheorem{lma}{Lemma}[section]

\newcommand{\beqa}{\begin{eqnarray}}
\newcommand{\eeqa}{\end{eqnarray}}

\newcommand{\pf}{\noindent {\bf Proof:} $\s$ }
\newcommand{\epf}{ \hfill$\diamondsuit$ \medskip}

\newcommand{\ds}{\displaystyle}
\newcommand{\beq}{\begin{equation}}
\newcommand{\eeq}{\end{equation}}
\newcommand{\lbl}{\label}
\newcommand{\s}{\; \;}

\newcommand{\la}{\lambda}

\newcommand{\ra}{\rightarrow}
\newcommand{\al}{\alpha}

\newcommand{\p}{\varphi}
\title[Infinitely many solutions]{Infinitely many solutions and asymptotics for resonant oscillatory problems }

\author[P.~Korman \& D.~Schmidt]{
Philip Korman   \\ 
Department of Mathematical Sciences \\ 
University of Cincinnati \\ 
Cincinnati Ohio 45221-0025 \\
\\
Dieter S. Schmidt \\
Department of Computer Science \\
University of Cincinnati \\
Cincinnati, Ohio 45221-0030
}

\dedicatory{Dedicated to the memory of Alan C. Lazer, a great mathematician and a dear friend}

\date{}

\begin{document}

\maketitle
\begin{abstract} 
For a class of oscillatory resonant problems, involving Dirichlet problems for semilinear PDE's on balls and rectangles in $R^n$, we show the existence of infinitely many solutions, and study the global solution set. The first harmonic of the right hand side is not required  to be zero, or small. We also derive   asymptotic formulas in terms of the first harmonic of solutions,  and illustrate their accuracy by numerical computations. The  numerical method is explained in detail.
 \end{abstract}

\begin{flushleft}
\end{flushleft}

\begin{flushleft}
\end{flushleft}

\section{Introduction}
\setcounter{equation}{0}
\setcounter{thm}{0}
\setcounter{lma}{0}

We study multiplicity of solutions for semilinear equations with linear part at resonance, a direction of research initiated by the classical paper of E.M. Landesman and A.C. Lazer \cite{L}. Several classes of oscillatory resonant problems on balls and rectangles in $R^n$ are considered. Our focus is on the existence of infinitely many solutions, which we establish by studying global solution  curves. Our results are supported by asymptotic analysis and numerical computations, and they extend related research in \cite{C},\cite{D1},\cite{SS},\cite{SS1}.
\medskip

Next we describe one of our main results, and the approach used. Let $B$ be the unit ball in $R^2$, $x^2+y^2<1$. For a class of oscillatory  resonant problems, with $h(u)=\sqrt{u} \sin \left[ \ln \left(u^{\frac{3}{2}} +1\right)\right]$ and $r=\sqrt{x^2+y^2}$,
\beqa
\lbl{0}
& \Delta u+\la _1 u+h(u)=g(x,y)=\mu _1 \p _1(r)+e(x,y) \s \mbox{for $x \in B$} \,, \\ \nonumber
& u=0 \s \mbox{on $\partial B$}
\eeqa
we show the existence of infinitely many solutions for any $g(x,y) \in L^2(B) \cap  C^{\al }(B)$, $\al >0$.  Here $ \Delta u=u_{xx}(x,y)+u_{yy}(x,y)$, while $(\la _1, \p _1(r))$ is the principal eigenpair of the Laplacian on $B$, with zero boundary conditions, $\mu _1 \in R$, $e(x,y) \in \p_1 ^{\perp}$ in $L^2(B)$, and $e(x,y) \in L^2(B) \cap C^{\al }(B)$, for some $\al \in (0,1)$. The function $e(x,y)$ is not assumed to be radially symmetric. Decompose solutions of (\ref{0}) as $u(x,y)=\xi _1 \p _1(r)+U(x,y)$, with $\xi _1 \in R$ and $U(x,y) \in \p_1 ^{\perp}$ in $L^2(B)$ (here $\p _1 ^{\perp}$ denotes the orthogonal complement of $\p _1$ in $L^2(B)$; similar decomposition is used throughout the paper).  We prove that the solution set of (\ref{0}) is exhausted by a continuous solution curve $(u(x,y), \mu _1)(\xi _1)$ parameterized by $\xi _1 \in R$, and  a section of this curve $\mu _1=\mu _1(\xi _1)$ oscillates toward $\pm \infty$ as $\xi _1 \ra \infty$, see Figure \ref{fig:0} below. We find this result to be rather surprising. A more typical behavior is that $\mu _1(\xi _1) \ra 0$  as $\xi _1 \ra \infty$. The choice of nonlinearity $h(u)$ was dictated by rather restrictive conditions that we needed to impose to obtain a continuous curve of solutions. For more natural nonlinearities we can only assert the existence of a continuum of solutions (as in the case of (\ref{2}) below), and in some cases we have to rely on the numerically produced solutions curves (e.g., for (\ref{1}) below), for which we can still derive accurate asymptotic formulas. In dimension $n=1$ a similar result was proved in \cite{K10} for $h(u)=u^p \sin u$, $\frac{1}{2}<p<1$. We remark that  solution curves of the form $(u(x,y), \mu _1)(\xi _1)$ appeared previously in R. Schaaf and K. Schmitt \cite{SS1}. 
\medskip

For a model resonant problem
\beqa
\lbl{1}
& \Delta u+\la _1 u+u \sin u=\mu _1 \p _1(r)+e(x,y) \s \mbox{for $x \in B$} \,, \\ \nonumber
& u=0 \s \mbox{on $\partial B$}
\eeqa
we provide asymptotical and computational evidence, as well as heuristic justification of the following conjecture: there exist two numbers $0<a<A$ so that the problem (\ref{1}) has infinitely many solutions for $\mu _1 \in (-a,a)$, there are at most finitely many solutions for $\mu _1$ outside of $(-a,a)$, and no solutions exist for $|\mu _1|>A$. A similar situation occurs for rectangles in two dimensions. The existence of infinitely many solutions for an interval of $\mu _1$'s (bounded or unbounded) is a new, and actually a rare phenomenon (as evidenced by the results of this paper, including  numerical computations). For the problem (\ref{1}) on balls and rectangles in dimensions  higher than two, and for ``most" other nonlinear terms,  infinitely many solutions occur only at $\mu _1=0$.
The restriction to the case $\mu _1=0$ is common in the literature, see e.g., \cite{C}, \cite{FN}, \cite{I}, \cite{SS}, \cite{SS1}. The nonlinear term in (\ref{1}) occurred previously in \cite{K17}, in connection with the oscillatory bifurcation from infinity.
\medskip

The problem (\ref{1}) can be seen as a limiting case of another model problem
\beqa
\lbl{2}
& \Delta u+\la _1 u+u^p \sin u=\mu _1 \p _1(r)+e(x,y) \s \mbox{for $x \in B$} \,, \\ \nonumber
& u=0 \s \mbox{on $\partial B$} \,,
\eeqa
with $p \in (0,1)$, to which the well known results of D. Costa et al \cite{C} apply at $\mu _1=0$ (the results of \cite{C} hold for more general domains), see also \cite{SS}. It follows from \cite{C} that for $\mu _1=0$  the problem (\ref{2}) has infinitely many solutions, and moreover $\ds \frac{u}{\max _Bu} \ra \p _1(r)$ for large solutions. We derive a rather precise  asymptotic formula for   $\mu _1=\mu _1(\xi _1)$ in case $|\xi _1|$ is large, and this formula tends to be accurate for small $|\xi _1|$ as well. The nonlinear term $u \sin u$ in (\ref{1}) has linear growth at infinity, while (\ref{2}) requires sublinear growth and is related to the result  in \cite{C}.
\medskip 

In addition to balls in $R^2$, we obtain asymptotic formulas and perform computations of solution curves for rectangular domains, and for radial solutions on balls  in any dimension.

\section{Global solution set for a ball in $R^2$}
\setcounter{equation}{0}
\setcounter{thm}{0}
\setcounter{lma}{0}

Let $J_0(z)$ be the Bessel function of order zero, with  $J_0(0)=1$, and denote by $\nu _1>0$ its first root, $\nu _1 \approx 2.405$. The principal eigenpair of the Laplacian on the unit ball $B \in R^2$ is $\la _1=\nu _1^2  \approx 5.78$, $\p _1(r)=c_0J_0(\nu _1 r)$ with $r=\sqrt{x^2+y^2}$, with $c_0$ chosen so that 
\[
\int _B \p _1^2(r) \, dxdy=2 \pi c_0^2 \int _0^1 J_0^2(\nu _1 r) \, r \, dr=1 \,, 
\]
which is
\beq
\lbl{2in}
c_0=\frac{1}{\sqrt{2 \pi \int _0^1 J_0^2(\nu _1 r) \, r\, dr}} \approx 1.09 \,.
\eeq 
Observe that $\p _1(0)=c_0$. We shall also need the second eigenvalue $\la _2$. Recall (see e.g. \cite{pr}) that the  eigenvalues of the Laplacian on $B$ with zero boundary condition are $\la _{n,m}=\al _{n,m}^2$ ($n=0,1,2,\ldots$; $m=1,2,\ldots$) with the corresponding eigenfunctions $J_n \left( \al _{n,m} r \right) \left(\al  \cos n \theta+\beta \sin  n \theta  \right)$, where $\al _{n,m}$ is the $m$-th root of $J_n(x)$, the $n$-th Bessel function ($\al$ and $\beta$ are arbitrary constants). One calculates $\la _2=\al _{1,1}^2  \approx 14.62$, with $\al _{1,1} \approx 3.83$, and $\p _2=J_1 \left( \al _{1,1} r \right) \left(\al  \cos  \theta+\beta \sin   \theta  \right)$. The principal  eigenvalue is simple, while all other eigenvalues have multiplicity two, because any two Bessel functions with indices differing by an integer do not have any roots in common, see G.N. Watson \cite{w}.
\medskip

Let us now recall the following result from \cite{K9} and \cite{K10}. It deals with PDE's on a general domain $ \Omega \subset R^n$
\beq
\lbl{intro1}
\s\s \Delta u+h(u)=\mu _1 \p _1(r)+e(x)  \s \mbox{for $x \in \Omega$}, \s u=0 \s \mbox{on $\partial \Omega$} \,.
\eeq
Here $x \in R^n$, $r=|x|$ and $(\la _1, \p _1(r))$ is the principal eigenpair of the Laplacian on $\Omega$, with zero boundary conditions, $\mu _1 \in R$, $e(x) \in \p_1 ^{\perp}$ in $L^2(\Omega)$, and $e(x) \in C^{\al }(\Omega) \cap L^2(\Omega)$, for some $\al \in (0,1)$.
\begin{thm}\lbl{thm:old}
For the problem (\ref{intro1}) assume that $h(u) \in C^2(R)$, $f(x) \in L^2(\Omega)$, and 
\beq
\lbl{2a1} 
h'(u) <\la _2- \la _1 , \s \mbox{for all $u \in R$}\,,
\eeq
\beq
\lbl{2a2}
|h(u)|<\gamma |u|+c , \s \mbox{with $0<\gamma <\la _2- \la _1$, $c \geq 0$, and  $u \in R$} \,.
\eeq
Then the solution set of (\ref{intro1}) consists of a continuous  curve  $(u(x), \mu _1)(\xi _1)$ parameterized by $\xi _1 \in R$.
\end{thm}

The following lemma extends a similar result in  \cite{K10}. Recall that  solutions of  (\ref{intro1}) are decomposed as $u(x)=\xi _1 \p _1(r)+U(x)$.
\begin{lma}\lbl{lma:old}
In the conditions of Theorem \ref{thm:old} assume that $\lim _{|u| \ra \infty} \frac{h(uz)}{u}=0$, uniformly in $z \in R$.  Then as $\xi _1 \ra \pm \infty$, the solutions of (\ref{intro1}) satisfy $\frac{u(x)}{\xi _1} \ra \p _1(x)$ in $H^1(\Omega)$ (and also in $L^2(\Omega)$). Moreover, if $h(\xi  _1 u)=O\left( |\xi  _1|^p \right)$ as $|\xi _1| \ra \infty$ uniformly in $u \in R$, then $||U(x)||_{H^1(\Omega)}=O\left( |\xi _1|^{\frac{p}{2}} \right)$ as $|\xi _1| \ra \infty$.
\end{lma}

\pf
By  Theorem \ref{thm:old} we have a solution curve $(u(x),\mu _1)(\xi _1)$. From (\ref{0})
\[
\Delta U+\la _1 U+h(\xi _1 \p _1+U)=\mu _1 \p _1+e \,.
\]
Letting here $U=\xi _1 V$, obtain
\[
\Delta V+\la _1 V=-\frac{h(\xi _1 \left(\p _1+ V \right))}{\xi _1}+\frac{\mu_1}{\xi _1} \p _1+ \frac{e}{\xi _1} \,.
\]
Multiplying by $V$ and integrating, we conclude that $\int _\Omega  |\nabla V|^2 \, dx=O\left( |\xi _1|^{p-1} \right)$, as $\xi _1 \ra \pm \infty$, and the lemma follows (since $V \perp \p _1$ in $L^2(\Omega)$, obtain $-\int_\Omega \Delta V V \, dx =\int _\Omega  |\nabla V|^2 \, dx \geq \la _2 \int _\Omega V^2 \, dx$).
\epf

\noindent Next we present one of our main results.

\begin{thm}
There exist $h(u) \in C^2(R)$ for which the problem (\ref{0}) has infinitely many positive solutions or any $g(x,y) \in L^2(B)\, \cap \, C^{\alpha }(B)$. Moreover, all solutions of  (\ref{0}) lie on a continuous solution curve $(u(x,y),\mu _1)(\xi _1)$, and $\mu _1(\xi _1)$ oscillates toward $\pm \infty$, as $\xi _1 \ra +\infty$.
\end{thm}
\pf 
We exhibit $h(u)$ satisfying the conditions of Theorem \ref{thm:old}, for which $\mu _1(\xi _1)$ oscillates toward $\pm \infty$, as $\xi _1 \ra +\infty$. Take $h(u)=\sqrt{u} \sin \left[ \ln \left(u^{\frac{3}{2}} +1\right)\right]$ for $u  \geq 0$ (observe that $h(u) \in C^2[0,\infty)$ and $h'(u)<1$ for $u \in [0,\infty)$), then  extend $h(u)$ arbitrarily to $(-\infty,0)$ so that $h(u) \in C^2(R)$ and the conditions (\ref{2a1}) and (\ref{2a2}) hold. Calculate (for $u \geq 0$) $h'(u)=\frac{\sin \left[ \ln \left(u^{\frac{3}{2}} +1\right)\right]}{2 \sqrt{u}}+\frac{3u\cos \left[ \ln \left(u^{\frac{3}{2}} +1\right)\right]}{2 \left(u^{\frac{3}{2}} +1\right)}$, and one of the anti-derivatives of $h(u)$:
\[
H(u)=\frac{1}{3} \left(u^{\frac{3}{2}} +1\right) \left(\sin \left[ \ln \left(u^{\frac{3}{2}} +1\right)\right]-\cos \left[ \ln \left(u^{\frac{3}{2}} +1\right)\right]     \right)
\]
\[
=\frac{\sqrt{2}}{3} \left(u^{\frac{3}{2}} +1\right) \sin \left[ \ln \left(u^{\frac{3}{2}} +1\right) -\frac{\pi}{4} \right] \,.
\]
Multiply the equation in (\ref{0}) by $\p _1$ and integrate over $B$:
\beq
\lbl{70}
\mu _1=\int _B h \left(\xi _1 \p _1+U(x,y) \right) \p _1 \,dxdy \,.
\eeq
Since $||U(x,y)||_{L^2(B)}=O(\xi _1^{\frac14})$ by Lemma \ref{lma:old}, and  $||h'(\xi _1 \p _1+U )||_{L^2(B)}=O(\frac{1}{\xi _1^{\frac12}})$ as $\xi _1 \ra \infty$, we have by the mean-value theorem
\[
||h \left(\xi _1 \p _1+U \right) - h \left(\xi _1 \p _1 \right)||_{L^2(B)}=o(1) \,.
\]
Then (\ref{70}) becomes
\beq
\lbl{70.2}
\mu _1 (\xi _1)=2 \pi \int _0^1 h \left(\xi _1 \p _1 \right)\p _1 r \,dr+o(1) \,.
\eeq
Denoting $f(r)=\frac{r\p _1(r)}{\p' _1(r)}$, and integrating by parts (using $\p _1(0)=c_0$)
\beq
\lbl{71}
\int _0^1 h \left(\xi _1 \p _1 \right)\p _1 r \,dr=\frac{1}{\xi _1} \int _0^1 f(r) \frac{d}{dr} H\left(\xi _1 \p _1 \right) \,dr
\eeq
\[
=-\frac{1}{\xi _1} \int _0^1 f'(r)  H\left(\xi _1 \p _1 \right) \,dr-\frac{1}{\xi _1} f(0)H\left(c_0\xi _1 \right) \,.
\]

The second term on the right  oscillates toward $\pm \infty$ as $\xi _1 \ra \infty$ with the amplitude approaching $\frac{\sqrt{2}}{3}c_0^{\frac{3}{2}}|f(0)| \xi _1^{\frac{1}{2}}$. We show next that it dominates the first term on the right in (\ref{71}). A computation shows that $f'(r)>0$ for all $r \in (0,1)$, and hence the first term is estimated in absolute value by a quantity approaching $\frac{\sqrt{2}}{3}c_0^{\frac{3}{2}} \xi _1^{\frac{1}{2}} \int _0^1 f'(r) J_0^{\frac{3}{2}} (\nu _1 r) \, dr$. Calculate $f(0)=\frac{\p _1(0)}{\p'' _1(0)}=-\frac{2}{\nu _1^2}$, so that $|f(0)|=\frac{2}{\nu _1^2} \approx 0.34$. (From $\p _1''+\frac{1}{r} \p _1'+\nu _1^2 \p_1=0$, it follows that $2 \p ''_1(0)=- \nu _1^2 \p _1(0)$.) Another computation shows that $\int _0^1 f'(r) J_0^{\frac{3}{2}} (\nu _1 r) \, dr \approx 0.1$, so that the second term dominates in (\ref{71}). Since $g(x,y) \in C^{\alpha }(B)$ the convergence $\frac{u(x,y)}{\xi _1} \ra \p _1(r)$ is in $C^2(B)$ by the elliptic regularity. It follows that $u(x,y)>0$ for large $\xi _1$, and hence it satisfies (\ref{0}) with the original $h(u)$ (before the extension). We conclude that for any $g(x,y) \in L^2(B) \cap  C^{\alpha }(B)$ the problem (\ref{0}) has infinitely many positive solutions.
\epf

\begin{figure}
\begin{center}
\scalebox{0.9}{\includegraphics{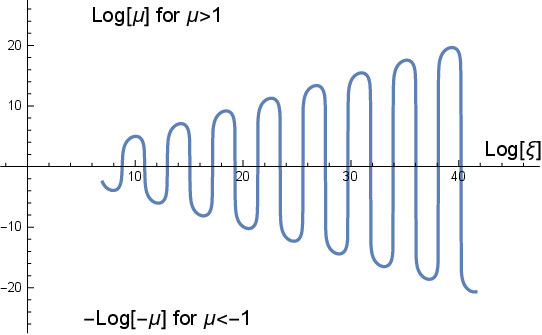}}
\end{center}
\caption{ The  solution curve $\mu _1= \mu _1(\xi _1)$ of (\ref{0}), oscillating to $\pm \infty$.  Values with $|\mu_1|<1$ are not shown.}
\lbl{fig:0}
\end{figure}

In Figure \ref{fig:0} we present an approximation of the solution curve $\mu _1=\mu _1(\xi _1)$ for the problem (\ref{0}), with $h(u)=\sqrt{u} \sin \left[ \ln \left(u^{\frac{3}{2}} +1\right)\right]$, computed using the formula (\ref{70.2}). We performed computations on very large intervals (along both $\xi _1$ and $\mu _1$ axes), and to make the resulting picture manageable a logarithmic scale is used for both $\xi _1$ and $\mu _1$. ($Log$ denotes the natural logarithm in {\em Mathematica}. Note that the curve in this presentation is not continuous.  It only appears so, since all values with $|\mu|<1$ have been omitted.) The result is the solution curve oscillating toward $\pm \infty$. In Figure \ref{fig:0a} we show the same curve in the original $(\xi _1, \mu _1)$ coordinates. It is not apparent from that picture that $\mu _1 (\xi _1) \ra \pm \infty$ as $\xi _1  \ra \infty$, since $\mu _1 (\xi _1)$ keeps the same sign on large, and ever increasing intervals.

\medskip

\begin{figure}
\begin{center}
\scalebox{0.9}{\includegraphics{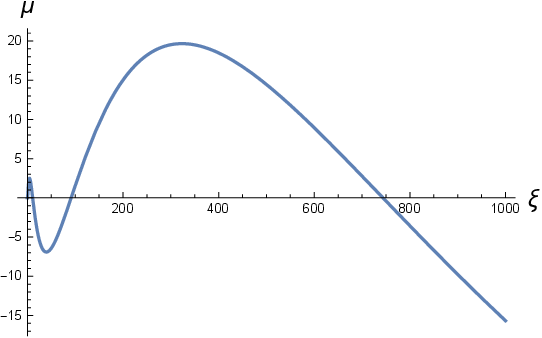}}
\end{center}
\caption{ The  solution curve $\mu _1= \mu _1(\xi _1)$ of (\ref{0}) for $\xi _1 \in (0,1000)$}
\lbl{fig:0a}
\end{figure}
\medskip

There are  other   $h(u)$  which can be handled similarly (including {\em Mathematica} being able to calculate the integral $H(u)$ in elementary functions). We mention $h(u)=u \sin \left[ \ln \left(u^{2} +1\right)\right]$ and $h(u)=\sin \left[ \ln \left(u +1\right)\right]$.
\medskip

Next we turn to more ``natural" nonlinearities $h(u)$.

\begin{thm}
If $0<p<1$ there is a continuum of solutions of (\ref{2}) $(u(x),\mu _1)(\xi _1)$ parameterized by the first harmonic $\xi _1 \in (-\infty,\infty)$. Along this continuum, $\lim _{\xi _1 \ra \pm \infty} \mu _1(\xi _1)=0$.  Moreover, the asymptotic formula (\ref{5}) below holds.
\end{thm}

\pf
The existence of a continuum of solutions of (\ref{2}) $(u(x),\mu _1)(\xi _1)$ parameterized by the first harmonic $\xi _1 \in (-\infty,\infty)$ follows by the result of R. Schaaf and K. Schmitt \cite{SS1}, which was based on E.N. Dancer \cite{D1}, see also P. Korman \cite{K10}. By the Lemma \ref{lma:old} and elliptic regularity it follows that $\frac{u}{\xi _1} \ra \p _1$ in $C^2(B)$, as $|\xi _1| \ra \infty$.
\medskip

We now derive an asymptotic formula for $\mu _1(\xi _1)$. Multiply the equation (\ref{2}) by $\p _1$, and integrate over $B$, then use Lemma \ref{lma:old} and elliptic regularity
\beqa
\lbl{3}
& \mu _1 =\int _B u^p \sin u \, \p _1  \, dx \, dy =\int_0^{2 \pi} \int _0^1  u^p \sin u \, \p _1 r \, dr \, d \theta \\
& =2 \pi\xi _1^p \int _0^1  \p _1^{p+1} \sin \left( \xi _1 \p _1 \right)   r \, dr+o(\xi _1^p)\,. \nonumber
\eeqa
Integration by parts gives
\beq
\lbl{4}
\int _0^1  \p _1^{p+1} \sin \left( \xi _1 \p _1 \right)   r \, dr= \int _0^1 \frac{\p _1 ^{p+1} r}{\xi _1 \p _1'} \frac{d}{dr} \left[ -\cos \left( \xi _1 \p _1 \right) \right] \, dr
\eeq
\[
=-\frac{1}{\xi _1} g(r) \cos \left(\xi _1 \p _1 \right) \, {\Huge |}_{_{0}}^{^{1}}+\frac{1}{\xi _1}\int _0^1 g'(r) \cos \left(\xi _1 \p _1 \right) \, dr \,,
\]
where $g(r)=\frac{\p _1^{p+1} \, r }{ \p '_1}$.  Observe that $g(1)=0$, while
\[
g(0)=c_0^{p+1}\lim _{r \ra 0} \frac{ r }{ \p '_1(r)}=\frac{ c_0^{p+1} }{ \p ''_1(0)}=-\frac{2c_0^p }{ \nu _1^2} \,.
\]
(From $\p _1''+\frac{1}{r} \p _1'+\nu _1^2 \p_1=0$, it follows that $2 \p ''_1(0)=-c_0 \nu _1^2$.)
Hence
\[
-g(r) \cos \left(\xi _1 \p _1 \right) \, {\Huge |}_{_{0}}^{^{1}}=-\frac{2 c_0^p \cos c_0 \xi _1}{\nu _1^2} \,.
\]
One checks that  $g'(r) \in C[0,1]$ is a bounded function. It follows that the   oscillating integral $\int _0^1 g'(r) \cos \left(\xi _1 \p _1 \right) \, dr$   is $o(1)$. From (\ref{3})
\beq
\lbl{5}
\mu _1 =-\frac{4 \pi \xi _1^{p-1} c_0^p \cos c_0 \xi _1 }{\nu _1^2} +o \left(\xi _1^{p-1} \right)\,.
\eeq
It follows that $\mu _1(\xi _1) \ra 0$ as $\xi _1 \ra \infty$, 
concluding the proof.
\epf
\smallskip

\noindent
{\bf Remark }
In limiting case $p=1$ the formula (\ref{5}) (while not rigorously justified) indicates   that $\mu_1(\xi _1)$ is  asymptotic to a multiple of $\cos c_0 \xi _1$, suggesting that there is a $\mu _0 >0$ so that for $|\mu |< \mu _0$ the problem  (\ref{2}) has infinitely many solutions. The same  conclusion is supported  by  our numerical computations, including  the following example.
\medskip

\noindent
{\bf Example} $\;$ We computed the solution curve $\mu _1= \mu _1(\xi _1)$ for the following example, with the linear part at resonance,
\beqa
\lbl{1ex}
& \Delta u+\la _1 u+u \sin u=\mu _1 \p _1(r)+xy \s \mbox{for $(x,y) \in B $} \,, \\ \nonumber
& u=0 \s \mbox{on $\partial B$} \,.
\eeqa
Observe that $\int _B xy \, \p_1(r) \, dxdy=0$.
The solution curve $\mu _1= \mu _1(\xi _1)$ (solid line) is presented in Figure \ref{fig:1}.  Notice a close agreement with the asymptotic formula (\ref{5}) (dashed line).

\begin{figure}
\begin{center}
\scalebox{0.95}{\includegraphics{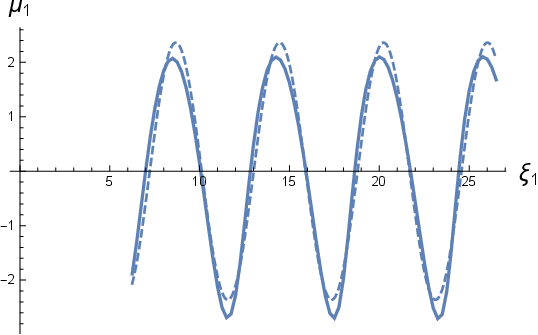}}
\end{center}
\caption{ The  solution curve $\mu _1= \mu _1(\xi _1)$ of (\ref{1ex}), compared with (\ref{5})}
\lbl{fig:1}
\end{figure}
\medskip

\noindent
{\bf Example} $\;$ We computed the solution curve $\mu _1= \mu _1(\xi _1)$ for an  example of (\ref{2}), with $p=\frac{1}{2}$,
\beqa
\lbl{1exa}
& \Delta u+\la _1 u+u^{\frac{1}{2}} \sin u=\mu _1 \p _1(r)+x^2y-3xy^4 \s \mbox{for $x \in B$} \,, \\ \nonumber
& u=0 \s \mbox{on $\partial B$} \,.
\eeqa
Observe that $\int _B \left(x^2y-3xy^4 \right) \, \p_1(r) \, dxdy=0$.
The solution curve $\mu _1= \mu _1(\xi _1)$ (solid line) is presented in Figure \ref{fig:2}. The solutions $u(x,y)$ {\em are not radially symmetric}, although {\em they get arbitrarily close to radially symmetric functions as $|\xi _1| \ra \infty$}, according to our results. Again, computations show a good  agreement with the asymptotic formula (\ref{5}) (dashed line).

\begin{figure}
\begin{center}
\scalebox{0.95}{\includegraphics{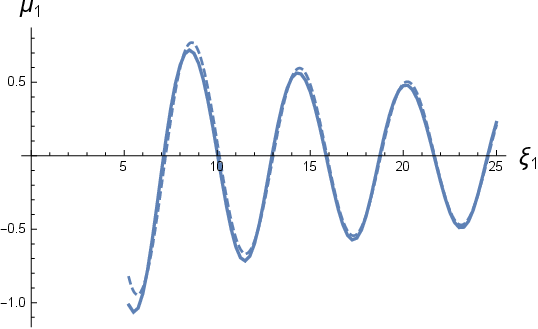}}
\end{center}
\caption{ The  solution curve $\mu _1= \mu _1(\xi _1)$ of (\ref{1exa}), compared with (\ref{5})}
\lbl{fig:2}
\end{figure}

\section{Asymptotic formula  in case of a rectangle}
\setcounter{equation}{0}
\setcounter{thm}{0}
\setcounter{lma}{0}

Let $R=\{0<x<a \} \times \{0<y<b \}$ be a rectangle in $R^2$. In this section we present computation of the solution curve $\mu _1=\mu _1(\xi  _1)$ on $R$ for the  nonlinearity considered above
\beqa
\lbl{r1}
& \Delta u+\la _1 u+u \sin u=\mu _1 \p _1(x,y)+e(x,y) \s \mbox{for $(x,y) \in R$} \,, \\ \nonumber
& u=0 \s \mbox{on $\partial R$} \,,
\eeqa
and derive an asymptotic formula for $\mu _1(\xi  _1)$. 
Here the principal eigenfunction $\p _1(x,y)=\frac{2}{\sqrt{ab}} \sin \frac{\pi}{a} x \sin \frac{\pi}{b} y$ satisfies $\int _R \p _1^2 \, dxdy=1$, and the corresponding  principal eigenvalue of $-\Delta$ is $\la _1=\frac{\pi^2}{a^2}+\frac{\pi^2}{b^2}$. It is assumed that  $\int \!\! \int _{R} e(x,y) \p _1 (x,y)\, dxdy=0$.
We decompose the solution of (\ref{r1}) as $u(x,y)=\xi _1 \p _1(x,y)+U(x,y)$, with $\int \!\! \int _{R} U(x,y) \p _1(x,y) \, dxdy=0$, and $\xi _1=\int \!\! \int _{R} u(x,y) \p _1(x,y) \, dxdy$.
\medskip

Unlike the case of a ball, we shall use the stationary phase method to derive the following asymptotic formula for $\xi _1$ large:
\beq
\lbl{r2}
\mu _1(\xi  _1) \sim \frac{4 \sqrt{ab}}{\pi} \sin \left(\frac{2}{\sqrt{ab}} \, \xi _1-\frac{\pi}{2} \right) \,,
\eeq
after we  recall the following known lemmas.

\begin{lma}\lbl{lma:as}
Assume that $f(x) \in C^2(x_0-a_0,x_0+a_0)$ for some constants $\al \,, a_0 >0$. Then as $\mu \ra \infty$
\[
\int_{x_0-a_0}^{x_0+a_0} f(x)e^{-i\al \mu (x-x_0)^2} \, dx=e^{-i \frac{\pi}{4}} \sqrt{ \frac{\pi}{\al \mu }}f(x_0)+O \left( \frac{1}{\mu} \right) \,.
\]
\end{lma}
This lemma follows from the part (i) of the following more general lemma, which we will also need, see e.g., p. 83 in \cite{K2} for the proof.
\begin{lma}\lbl{lma:new}
(i) Assume that  $f(x)$ and $g(x)$ are of class $C^2[a,b]$ and $g(x)$ has a unique critical point $x_0$ on $[a,b]$, and moreover $x_0 \in (a,b)$ and $g''(x_0) \ne 0$ (so that $x_0$ gives a global max or global min on $[a,b]$). Then as $\mu \ra \infty $ the following asymptotic formula holds
\[
\int _a^b f(x) e^{i \mu g(x)} \, dx=e^{i \left[ \mu g(x_0) \pm \frac{\pi}{4} \right]} \sqrt{ \frac{2 \pi}{\mu |g''(x_0)|}} 
\, f(x_0)+O \left(\frac{1}{\mu} \right) \,, 
\]
where one takes ``plus" if $g''(x_0)>0$ and ``minus"  if $g''(x_0)<0$.
\medskip

\noindent
(ii) Assume that the functions  $f(x)$ and  $g(x)>0$ are of class  $C^2[0,1]$, and satisfy
\[
g'(x)<0 \quad {\rm for \, all}\quad x\in(0,1], \quad {\rm and }\quad g'(0)=0,\,\,g''(0)<0\,.
\]
Then, as $\mu\to \infty$,
\[
\int_0^1 f(x)e^{i\mu g(x)}dx =  e^{i(\mu g(0)-
\frac{\pi}{4})}\sqrt {\frac{\pi}{2\mu |g''(0)|}}f(0) + O\left(\frac{1}{\mu}\right)\,.
\]
\end{lma}

Turning to the derivation of (\ref{r2}),
multiply (\ref{r1}) by $\p _1$ and integrate over $R$. Then use that $u \sim \xi _1 \p _1$  for $\xi _1$ large ($\sin u \sim \sin  \xi _1 \p _1$ by the mean-value theorem)
\[
\mu _1=\int \!\! \int _{R} u \sin u \, \p _1  \, dxdy \sim \xi _1 \int \!\! \int _{R}  \p _1 ^2 \sin  \xi _1 \p _1 \, dxdy=\xi _1 {\rm Im }\int \!\! \int _{R}  \p _1 ^2 e^{i \xi _1 \p _1}  \, dxdy
\]
\[
=\frac{4 \xi _1}{ab} {\rm Im }\int \!\! \int _{R} \sin ^2  \frac{\pi}{a} x  \sin ^2  \frac{\pi}{b} y \, e^{i \frac{2}{\sqrt{ab}}\, \xi _1 \sin \frac{\pi}{a} x \sin \frac{\pi}{b} y}  \, dxdy \,.
\]
For $\xi _1$ large there are fast oscillations around zero, except near the stationary point $(x_0,y_0)=(\frac{a}{2},\frac{b}{2})$. The approximation of the integral near the stationary point provides the dominant contribution to this integral. Using Taylor's formula near $(x_0,y_0)$
\[
\sin   \frac{\pi}{a} x  \sin  \frac{\pi}{b} y \approx 1-\frac{\pi ^2}{2 a^2}(x-x_0)^2-\frac{\pi ^2}{2 b^2}(y-y_0)^2 \,,
\]
because the second mixed partials vanish at $(x_0,y_0)$.  Then
\[
\mu _1 \sim \frac{4 \xi _1}{ab} {\rm Im } \left( e^{i \frac{2}{\sqrt{ab}}\, \xi _1} \int_0^a \sin ^2  \frac{\pi}{a} x \,  e^{-i \xi _1 \al _1(x-x_0)^2} \, dx   
\int_0^b \sin ^2  \frac{\pi}{b} y \, e^{-i \xi _1 \al _2(y-y_0)^2} \, dy 
\right),
\]
where $\al _1=\frac{\pi ^2}{ a^2\sqrt{ab}}$ and $\al _2=\frac{\pi ^2}{ b^2\sqrt{ab}}$. Using Lemma \ref{lma:as}, for large $\xi _1$,
\[
\mu _1 \sim \frac{4 \xi _1}{ab} {\rm Im } \left( \, e^{i \frac{2}{\sqrt{ab}}\, \xi _1} e^{-i \frac{\pi}{4}} \sqrt{ \frac{\pi}{\al _1  \xi _1 }}\, e^{-i \frac{\pi}{4}} \sqrt{ \frac{\pi}{\al _2 \xi _1 }}
\right) 
=\frac{4 \sqrt{ab}}{\pi} \sin \left(\frac{2}{\sqrt{ab}} \, \xi _1-\frac{\pi}{2} \right)\,.
\]

This formula, as well as our numerical calculations, suggests that  there exist two constants $0<a<A$ so that the problem (\ref{r1}) has infinitely many solutions for $\mu _1 \in (-a,a)$, there are at most finitely many solutions for $\mu _1$ outside of $(-a,a)$, and no solutions exist for $|\mu _1|>A$.
\medskip

\noindent
{\bf Example } On the rectangle $R_1=\{0<x<1 \} \times \{0<y<2 \}$ we computed the solution curve $\mu _1=\mu _1(\xi  _1)$ for the problem 
\beqa
\lbl{r3}
& \s \Delta u+\la _1 u+u \sin u=\mu _1 \p _1(x,y)+(x-\frac{1}{2})(y-1) \s \mbox{for $(x,y) \in R_1$} \,, \\ \nonumber
& u=0 \s \mbox{on $\partial R_1$} \,,
\eeqa
with $\la _1=\frac{5\pi^2}{4}$ and $\p _1(x,y)=\sqrt{2} \sin \pi x \sin \frac{\pi}{2} y$. Observe that
$\int \!\! \int _{R} (x-\frac{1}{2})(y-1) \p _1(x,y) \, dxdy=0$.
The solution curve $\mu _1= \mu _1(\xi _1)$ (solid line) is presented in Figure \ref{fig:3}. Once again, we have  a close agreement with the asymptotic formula (\ref{r2}) (dashed line).
\medskip

As in the case of balls, oscillations of $\mu _1(\xi _1)$ are decaying in the dimensions $n>2$, as will follow from the asymptotic formula that we present next. Consider the  $n$-dimensional rectangle $R=(0,a_1) \times (0,a_2)\times \cdots \times (0,a_n)$, and the problem (\ref{lma:as}) on $R$. The principal eigenfunction of the Laplacian on $R$ with Dirichlet boundary conditions, and satisfying $\int _{R} \p _1^2 \, dx=1$, is 
\[
\p _1=\frac{2^{\frac{n}{2}}}{\sqrt{a_1a_2 \cdots a_n}} \, \sin \frac{\pi}{a_1} x_1  \, \sin \frac{\pi}{a_2}  x_2 \cdots \, \sin \frac{\pi}{a_n}  x_n \,.
\]
As in the above derivation, one shows that 
\[
\mu _1 (\xi _1) \sim \frac{2^{\frac{n}{2} (3-\frac{n}{2})} \left(a_1a_2 \cdots a_n \right)^{\frac{n}{4}}}{\pi ^{\frac{n}{2}}}\xi _1^{1-\frac{n}{2}} \sin \left(\frac{2^{\frac{n}{2}}}{\sqrt{a_1a_2 \cdots a_n}} \, \xi _1-n\frac{\pi}{4} \right) \ra 0 \,,
\]
as $\xi _1 \ra \infty$, for $n>2$.
\begin{figure}
\begin{center}
\scalebox{0.95}{\includegraphics{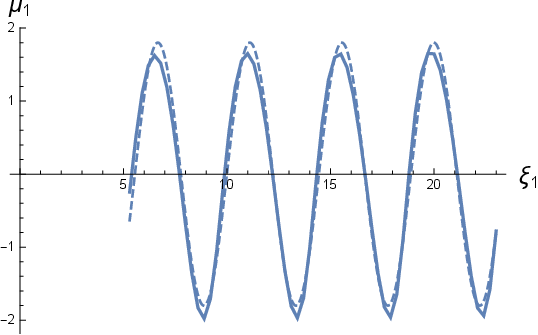}}
\end{center}
\caption{ The  solution curve $\mu _1= \mu _1(\xi _1)$ of (\ref{r3}), compared with (\ref{r2})}
\lbl{fig:3}
\end{figure}

\section{Numerical computations}
\setcounter{equation}{0}
\setcounter{thm}{0}
\setcounter{lma}{0}

We now describe the {\em Mathematica} program that was used to produce the solution curves $\mu _1=\mu _1(\xi _1)$, presented above. To avoid ambiguity of notation, in this section we shall write  $\xi$ instead of $\xi _1$ and $\mu$ instead of $\mu _1$, so that $u=\xi \p _1+U$, and $\mu=\mu (\xi)$. 
\medskip

Our program  handles semilinear equations on rather  general  domains solving
\beqa
\lbl{n1}
& \Delta u+f(u)=\mu  \p _1(x,y)+e(x,y) \s \mbox{for $(x,y) \in \Omega$} \,, \\ \nonumber
& u=0 \s \mbox{on $\partial \Omega$}
\eeqa
on bounded domains $\Omega \in R^2$, including rectangles and ellipses in $R^2$, and radially symmetric solutions on balls  in $R^n$. Here $ \Delta u=u_{xx}(x,y)+u_{yy}(x,y)$, while $(\la _1, \p _1)$ is the principal eigenpair of the Laplacian on $\Omega$, with zero boundary conditions and $\int _{\Omega} \p _1^2 \, dxdy=1$, $\mu  \in R$, $e(x,y) \in \p_1 ^{\perp}$ in $L^2(\Omega)$. 
Choosing a step size $h$ and an initial value $\xi =\xi _0$,  let $\xi _n=\xi _0+nh$. We are looking for $\mu =\mu _n$ for which the problem (\ref{n1}) has a solution $u(x) \equiv u_n(x,y)$, with 
\beq
\lbl{n2}
\int _{\Omega} u(x,y) \p _1(x,y)\, dxdy=\xi _n \,.
\eeq
We utilize {\em Mathematica}'s ability (the NDSolve command) to solve linear Dirichlet problems of the type
\[
\Delta u+a(x,y)u=b(x,y) \s \mbox{for $(x,y) \in \Omega$} \,, \s u=0 \s \mbox{on $\partial \Omega$} 
\]
on some domains $\Omega \in R^2$, including ellipses around the origin, and rectangles.
\medskip

Assuming that $(\mu _n,u_n (x,y))(\xi _n)$ is already computed, we use Newton's method to calculate  $(\mu _{n+1},u_{n+1} (x,y))(\xi _{n+1})$. We  calculate $(\mu _{n+1},u_{n+1} (x,y))$ using a sequence of iterates $(\mu ^k,u^k(x,y))$ beginning with $(\mu ^0,u^0(x,y))=(\mu _n,u_n (x,y))$. Assuming that $(\mu ^k,u^k (x,y))$ is already computed, we approximate $f(u) \approx f(u^k)+f'(u^k)(u-u^k)$, and solve the linear problem
\beqa \nonumber
& \Delta w+f'(u^k)w=\mu _1 \p _1(x,y)+ f'(u^k)u^k-f(u^k)+e(x,y) \s \mbox{for $(x,y) \in \Omega$} \,, \\ \nonumber
& w=0 \s \mbox{on $\partial \Omega$}  \nonumber
\eeqa
for $\mu$ and $w$ by the algorithm that is described next. One can decompose the solution in the form $w(x,y)=\mu w_1(x,y)+w_2(x,y)$, where $w_1$ and $w_2$ are solutions of  
\[
\Delta w_1+f'(u^k)w_1= \p _1 (x,y) \,, \s \mbox{for $(x,y) \in \Omega$} \,, \s w_1=0 \s \mbox{on $\partial \Omega$} \,,
\]
\beqa 
\lbl{n2.1}
& \Delta w_2+f'(u^k)w_2=  f'(u^k)u^k-f(u^k)+e(x,y) \,, \s \mbox{for $(x,y) \in \Omega$} \,, \\ \nonumber
& w_2=0 \s \mbox{on $\partial \Omega$} \,. \nonumber
\eeqa
After calculating $w_1$ and $w_2$, we look for $\mu$ such that $\int _{\Omega} w(x,y) \p _1(x,y)\, dxdy=\xi _{n+1}$, and declare that value of $\mu$ to be our new iterate $\mu ^{k+1}$, so that
\beq
\lbl{n5}
\mu ^{k+1}=\frac{\xi _{n+1}-\int _{\Omega} w_2(x,y) \p _1(x,y)\, dxdy}{\int _{\Omega} w_1(x,y) \p _1(x,y)\, dxdy} \,.
\eeq
The corresponding $w(x,y)$ is our next iterate $u^{k+1}(x,y)=\mu  ^{k+1} w_1(x,y)+w_2(x,y)$. The iterations are stopped once the relative error $\frac{\mu  ^{k+1}-\mu  ^{k}}{\mu  ^{k}}$ is small enough.
\medskip

Controlling the accuracy of  iterates is the major improvement of the present algorithm, compared with the one we used in \cite{KS}.
\medskip

\noindent
{\bf Remark } There is a  better way of choosing the initial iterate $u^0(x,y)$ at $\xi _{n+1}$ (corresponding to the ``predictor" of the predictor-corrector method). Write $u=u(x,y,\xi)$. Approximate $u(x,y,\xi_{n+1}) \approx u(x,y,\xi_{n})+u_{\xi}h=u_n(x,y)+u_{\xi}h$. Differentiate the equation {\ref{n1}) in $\xi$, and compare the result with the first formula in (\ref{n2.1}) to get: $ u_{\xi}=\mu '(\xi _n)w_1 \approx \frac{\mu _{n+1}-\mu _n}{h}w_1$. So that we take $u^0(x,y)=u_n(x,y)+\left(\mu _{n+1}-\mu _n \right)w_1$, using the last  function $w_1$ computed at $\xi _n$. Our experiments show  considerably faster convergence.

\section{Asymptotics and numerics for radial solutions}
\setcounter{equation}{0}
\setcounter{thm}{0}
\setcounter{lma}{0}

Let $B$ denote the unit ball around the origin  in $R^n$, $x \in R^n$ and $r=|x|$. We present computations of the solution curve $\mu _1=\mu _1(\xi  _1)$  for the  radial solutions $u=u(r)$ of the model problem 
\beqa
\lbl{rad1}
&\s u''(r)+\frac{n-1}{r}u'(r)+\la _1 u+ \sin u=\mu _1 \p _1(r)+e(r) \,, \s \mbox{for $0<r<1$}  \\ \nonumber
& u'(0)=u(1)=0 \,,
\eeqa
and derive an asymptotic formula for $\mu _1(\xi  _1)$. By the Theorem \ref{thm:old} there exists a continuous solution curve $(u(r),\mu _1)(\xi _1)$ that exhausts the solution set of (\ref{rad1}), and moreover $\frac{u(r)}{\xi _1} \ra \p _1$ in $C^2(B)$ as  $\xi _1 \ra \infty$. (The same result holds if one replaces $e(r)$ by $e(x) \in \p _1^{\perp}$.) Restricting to radial solutions we can perform numerical computations for $n \geq 3$.) 
Recall that  the principal eigenfunction of the Laplacian on $B$ is  $\p _1(r)=c_0 r^{-\frac{n-2}{2}} J_{\frac{n-2}{2}}(\nu _1 r)$, where  $\nu _1$ denotes the first root of the Bessel function $J_{\frac{n-2}{2}}(r)$, and 
$c_0$ is chosen so that 
\[
\int _B \p _1^2(r) \, dx=\omega _n c_0^2 \int _0^1 J^2_{\frac{n-2}{2}}(\nu _1 r) \, r \, dr=1 \,, 
\]
which is
\[
c_0=\frac{1}{\sqrt{\omega _n \int _0^1 J^2_{\frac{n-2}{2}}(\nu _1 r) \, r\, dr}}  \,.
\]
Here $\omega _n=\frac{n \pi ^{\frac{n}{2}}}{\Gamma(\frac{n}{2}+1)}$ gives the area of the unit ball in $R^n$, where $\Gamma$ denotes the gamma function.
The corresponding  principal eigenvalue is $\la _1=\nu _1 ^2$. It is assumed that  $\int \!\! \int _{B} e(r) \p _1 (r) \, dx=\omega _n \int _0^1 e(r) \p _1 (r) r^{n-1} \, dr=0$.
We decompose the solution of (\ref{r1}) as $u(r)=\xi _1 \p _1(r)+U(r)$, with $\int \!\! \int _{B} U(r) \p _1(r) \, dx=0$, and $\xi _1=\int \!\! \int _{B} u(r) \p _1(r) \, dx$. 
\medskip

We now derive  an asymptotic formula for the function $ \mu _1(\xi _1)$, by using that $u(x) \sim \xi _1 \p _1$ for large $\xi _1$. Multiplying the PDE version of the equation (\ref{rad1}) by $\p _1$ and integrating over $B$
\beqa \nonumber
& \mu _1=  \omega _n \int _0^1 \left( \sin u(r) \right) \p _1(r) \, r^{n-1} \,dr \sim  \omega _n \int _0^1 \left( \sin \xi _1 \p _1 \right) \p _1(r) \, r^{n-1} \,dr \\ \nonumber
& =\frac{\omega _n}{\xi _1} \int _0^1 \frac{\p _1(r) \, r^{n-1}}{\p' _1(r)} \frac{d}{dr} \left(- \cos \xi _1 \p _1 \right) \, dr \,. \nonumber
\eeqa
Integrating by parts, and denoting $f_1(r) \equiv \frac{\p _1(r) \, r^{n-1}}{\p' _1(r)}$, obtain ($\p_1 (1)=0$)
\beq
\lbl{rad2}
\mu _1 \sim \frac{\omega _n}{\xi _1} \left( f_1(0) \cos \xi _1 \p _1 (0) +\int _0^1  f'_1(r)\cos \xi _1 \p _1(r) \,dr \right) \,.
\eeq
The next steps depend on the dimension $n$.
\medskip

Assume that $n=2$. Then $\omega _2=2\pi$, $\p _1(r)=c_0 J_0(\nu _1r )$, where $\nu _1$ is the first root of the Bessel's function $J_0(x)$, $c_0$ is given by (\ref{2in}),  
$f_1(r) = \frac{\p _1(r) \, r}{\p' _1(r)}$, and 
$
f_1(0) = \frac{\p _1(0)  }{\p'' _1(0)} =-\frac{2}{\nu _1^2}
$.
We conclude that  for large $\xi _1$  (using that $\p_1 (0)=c_0$, and that the oscillatory integral in (\ref{rad2}) is $o(1)$)
\beq
\lbl{rad4}
\mu _1 \sim -\frac{4 \pi}{\xi _1 \nu _1^2} \cos c_0 \xi _1 \,,
\eeq
which is consistent with (\ref{5}) at $p=0$.
\medskip

Assume now that $n=3$. Then $\omega _3=4\pi$. Since $J_{\frac{1}{2}}(x)=\frac{\sin x}{\sqrt{x}}$, we have $\nu _1=\pi ^2$, $\p _1(r)=\frac{1}{\sqrt{2 \pi}} \frac{\sin \pi r}{r}$, $\p _1(0)=\sqrt{\frac{\pi}{2}}$. Also, $f_1(r)=\frac{\p _1 r^2}{\p' _1}$. It follows that $f_1(0)=0$, so that the first term in (\ref{rad2}) vanishes. Then by Lemma \ref{lma:new}  (ii),
\beqa \nonumber
& \mu _1 \sim \frac{4 \pi}{\xi _1} {\rm Re} \, \int _0^1 f'_1 e^{i \xi _1 \p _1} \, dr= \frac{4 \pi}{\xi _1} {\rm Re} \, e^{i \left(\xi _1 \p _1(0)-\frac{\pi}{4} \right)} \sqrt{\frac{\pi}{2\xi _1 |\p''_1(0)|}} f'_1(0)\\ \nonumber
& =\frac{4 \pi}{\xi _1} \cos  \left(\xi _1 \sqrt{\frac{\pi}{2}}-\frac{\pi}{4} \right) \sqrt{\frac{\pi}{2\xi _1 |\p''_1(0)|}} f'_1(0) \,.\nonumber
\eeqa

As above, calculate $3\p''_1(0)=-\pi ^2 \p_1(0)$, so that $|\p''_1(0)|=\frac{\pi ^{\frac52}}{3 \sqrt{2}}$. A short calculation shows that
\[
f'_1(0) = \frac{\p _1(0)  }{\p'' _1(0)} =-\frac{3}{\pi ^2}\,.
\]
We conclude that for large $\xi _1$
\beq
\lbl{rad5}
\mu _1 \sim -\frac{12 \sqrt{3 \sqrt{2}} }{\sqrt{2} \xi _1 ^{\frac32} \pi^{\frac74}} \cos  \left(\xi _1 \sqrt{\frac{\pi}{2}}-\frac{\pi}{4} \right) \,.
\eeq

For numerical computations of radial solutions we used the algorithm described in a previous section. Computations can be performed accurately in any dimension, because the linearization of (\ref{rad1}) involves a boundary value problem for an ODE, readily handled by {\em Mathematica}. The formula (\ref{n5}) takes the form
\[
\mu^{k+1}=\frac{\xi _{n+1}-\omega _n \int _0^1 w_2(r)\p _1(r) \, r^{n-1} \, dr}{\omega _n \int _0^1 w_1(r)\p _1(r) \, r^{n-1} \, dr} \,.
\]

\noindent
{\bf Example } For $n=3$ we solved the problem (\ref{rad1}) with $e(r)=\frac{\cos \pi r}{r}$ (observe that $\int _B e(r) \p _1(r) \, dx =0$), see Figure \ref{fig:5}. The solution of (\ref{rad1}) (solid line) is in a good agreement with the asymptotic formula (\ref{rad5}) (dashed line).

\begin{figure}
\begin{center}
\scalebox{0.95}{\includegraphics{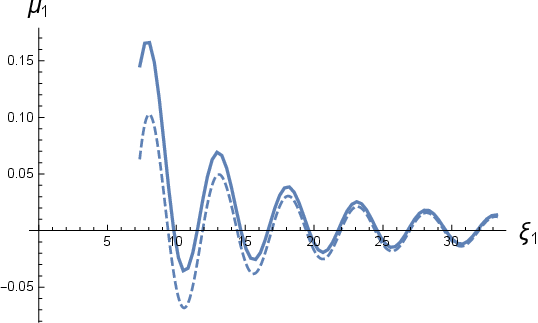}}
\end{center}
\caption{ The  solution curve $\mu _1= \mu _1(\xi _1)$ of (\ref{rad1}), $n=3$, compared with the asymptotic formula (\ref{rad5})}
\lbl{fig:5}
\end{figure}

\end{document}